\documentclass[oneside]{amsart}
\usepackage{amsfonts}
\usepackage{amssymb}
\usepackage{amsxtra}
\usepackage{amstext}
\usepackage[english]{babel}

\newcommand{\R}{\mathbb{R}}

\newcommand{\N}{\mathbb{N}}
\newcommand{\E}{\mathbb{E}}
\newcommand{\Z}{\mathbb{Z}}

\newcommand{\pp}{\mathbb{P}}

\newcommand{\kC}{\mathcal{C}}

\newcommand{\kO}{\mathcal{O}}

\newcommand{\kF}{\mathcal{F}}

\newtheorem {lem} {Lemma} [section]
\newtheorem {prop} {Proposition} [section]
\newtheorem {theo} {Theorem} [section]
\newtheorem {cor} {Corollary} [section]
\newtheorem {rem} {Remark} [section]

\newtheorem {Hypo} {Hypothesis} [section]

\title[On some generalized reinforced random walk on integers]
      {On some generalized reinforced random walk on integers}

\author{Olivier RAIMOND}
\address{Laboratoire Modal'X, Universit\'e Paris Ouest Nanterre La D\'efense, B\^atiment G, 200 avenue de la R\'epublique 92000 Nanterre, France.}
\email{olivier.raimond@u-paris10.fr}
\author{Bruno SCHAPIRA}
\address{D\'epartement de Math\'ematiques, B\^at. 425, Universit\'e Paris-Sud 11, F-91405 Orsay, cedex, France. }
\email{bruno.schapira@math.u-psud.fr}

\begin{document}

\begin{abstract} We consider Reinforced Random Walks where transitions probabilities are a function of the proportions of
times the walk has traversed an edge. We give conditions for
recurrence or transience. A phase transition is observed, similar to
Pemantle \cite{Pem000} on trees.
\end{abstract}

\maketitle


\bigskip

\noindent \textbf{A.M.S. classification:} 60F20.

\section{Introduction}
Let $G=(V,E)$ be a graph with $V$ the set of vertices and $E$ the
set of edges. The graph distance is denoted by $d$. Linearly reinforced
random walks $(X_n,n\ge 0)$ are nearest neighbor walks on $G$ (i.e. $X_n\in V$ and
$(X_n,X_{n+1})\in E$) whose laws are defined as follows: denote by
$(\mathcal{F}_n)_{n\ge 0}$ the natural filtration associated to
$(X_n,n\ge 0)$ and for $(x,y)\in E$, set
$$a_n(x,y)=a_0(x,y)+\Delta \sum_{i=1}^n 1_{\{(X_{i-1},X_i)=(x,y)\}},$$
with $a_0(x,y)>0$ and $\Delta>0$. Then for $(x,y)\in E$ and on the event $\{X_n=x\}$,
$$\pp(X_{n+1}=y\mid \kF_n)=\frac{a_n(x,y)}{\sum_{\{z \mid (x,z)\in E\}}a_n(x,z)}.$$
When $E$ is non-oriented (i.e. $(x,y)\in E$ implies $(y,x)\in E$),
for $(x,y)\in E$, set
$$b_n(x,y)=b_0(x,y)+\Delta \sum_{i=1}^n 1_{\{(X_{i-1},X_i)\in\{(x,y),(y,x)\}\}},$$
with $b_0(x,y)=b_0(y,x)>0$ and $\Delta>0$. The law of an undirected
reinforced random walk is such that for $(x,y)\in E$, on the event
$\{X_n=x\}$,
$$\pp(X_{n+1}=y\mid \kF_n)=\frac{b_n(x,y)}{\sum_{\{z \mid (x,z)\in E\}}b_n(x,z)}.$$

The law of a directed linearly reinforced random walk is the same as the law
of a random walk in a random environment. Indeed it is equivalent to
attach independent Polya urns to all sites and then De Finetti
Theorem implies that it is equivalent to attach independent random
probability vectors to each site. These probability vectors give the
transition probabilities for the walk (when located at this site).
When the graph is a non-oriented tree, the undirected reinforced
random walk, with initial position $\rho$, has the same law as a
directed reinforced random walk on $G$ with $a_0(x,y)=b_0(x,y)$ if
$d(\rho,y)=d(\rho,x)+1$ and $a_0(x,y)=b_0(x,y)+\Delta$ if
$d(\rho,y)=d(\rho,x)-1$ and reinforcement parameter $2\Delta$, or with
$a_0(x,y)=b_0(x,y)/(2\Delta)$ if $d(\rho,y)=d(\rho,x)+1$ and
$a_0(x,y)=(b_0(x,y)+\Delta)/(2\Delta)$ if $d(\rho,y)=d(\rho,x)-1$
and reinforcement parameter 1. This representation was first observed
by Coppersmith and Diaconis \cite{CD}.

For this class of models results on random walks in random
environment can be applied. When the graph is $\Z$, Solomon's
theorem shows that (directed and undirected) reinforced random walks
are a.s. recurrent when $a_0(x,x+1)=a_0(x,x-1)=a_0>0$ or
$b_0(x,x+1)=b_0>0$. When the graph is the binary tree and
$b_0(x,y)=b_0>0$, the undirected reinforced random walk is transient
for small $\Delta$ (or equivalently for large $b_0$) and recurrent
for large $\Delta$ (or equivalently for small $b_0$). This last
result was proved by R. Pemantle in \cite{Pem000}.

In this paper we address the question (posed by M. Bena\"{\i}m to
one of the authors) of what happens when the graph is $\Z$ and when on
the event $\{X_n=x\}$,
$$\pp(X_{n+1}=x+1\mid\kF_n)=f\left(\frac{a_n(x,x+1)}{a_n(x,x-1)+a_n(x,x+1)}\right),$$
where $f:[0,1]\to (0,1)$ is a smooth function. For general functions
$f$, these walks are no longer random walks in random environment.
So we have to use different techniques. But one can still attach to
each site independent urn processes (generalized Polya urns). Under
the assumption that the number of fixed points of $f$ is finite, if
the walk is recurrent, stochastic algorithms techniques show that
for all $x$, $a_n(x,x+1)/(a_n(x,x-1)+a_n(x,x+1))$ converges a.s.
toward a random variable $\alpha_x$. This random variable takes its
values in the set of fixed point of $f$. If
$a_0(x,x+1)=a_0(x,x-1)=a_0>0$, the sequence $(\alpha_x,x\in \Z)$ is
i.i.d. Let us remark that Solomon's theorem states that the random
walk in the random environment $(\alpha_x,x\in \Z)$ is a.s.
recurrent if and only if $\E[\ln(\alpha_x/(1-\alpha_x))]=0$.

We focus here on cases when either $f$ has a unique fixed point or all the fixed points are greater or
equal to $1/2$ and $f\geq 1/2$ on $[1/2,1]$. We
particularly study the case $a_0(x,x+1)=a_0(x,x-1)=a_0>0$. We give
criteria for recurrence and transience: 
\begin{itemize}
\item when there exists one fixed
point greater than $1/2$ the walk is transient and 
\item when $1/2$ is the
unique fixed point, depending on the initial condition $a_0$ and on
the shape of $f$ around $1/2$, the walk can be either recurrent or
transient. 
\end{itemize}
This last result shows that Solomon's criterion applied to the limiting values $(\alpha_x, x\in \Z)$ does not determine recurrence versus transience. The proofs of the theorems given here involve martingale
techniques inspired by the work of Zerner on multi-excited random
walks on integers \cite{Zer,Zer2}.

The paper is organized as follows: in section 2 reinforced random
walks are defined and their representation with urn processes is
given. In section 3 are given the results on urns that are needed to
prove the theorems given in sections 5 and 6. A zero-one law is
proved in section 4: recurrence occurs with probability 0 or 1. In
section 5 and 6 the case $f\geq 1/2$ and the case when there is a unique
fixed point are studied. The last section develops some examples.

\section{Notation}
\label{prelim}
\subsection{The urn model}
We consider an urn model where balls of the urn are only of two
types or colors, let say Red and Blue. Given a function $f:[0,1]\to
(0,1)$ we alter the draw by choosing at each step a Red ball with
probability $f(\alpha)$, if $\alpha$ is the proportion of Red balls.
Then we put back two Red (respectively Blue) balls in the urn if a
Red (respectively Blue) ball was drawn. In other words an urn
process associated to $f$ is a Markov process $((\alpha_n,l_n),n\ge
0)$  on $[0,1]\times (0,+\infty)$, where the transition
probabilities are defined as follows: for all $n$, $l_{n+1}=l_n+1$
and $\alpha_{n+1}$ is equal to $(l_n\alpha_n+1)/(l_n+1)$ with
probability $f(\alpha_n)$, or equal to $l_n\alpha_n/(l_n+1)$ with
probability $1-f(\alpha_n)$. In fact $\alpha_n$ represents the
proportion of Red balls and $l_n$ the total number of balls in the
urn at time $n$ (at least if $l_0$ and $\alpha_0 l_0$ are integers).
By abuse of notation we will sometime call the first coordinate
$(\alpha_n,n \ge 0)$ an urn process associated to $f$ (if there is no ambiguity on $l_0$). 
This model was introduced in \cite{HLS}, and then further studied in particular by Pemantle \cite{Pem1}, Duflo \cite{D} and Bena\"im and Hirsch (see \cite{BH} and \cite{B}). We refer also the reader to the survey \cite{Pem2} section $2.4$, $2.5$ and $3.2$ for more details and references.

\subsection{Generalized reinforced random walks}
Here we consider a particular model of (directed) reinforced random
walk $(X_n,n\ge 0)$ on $\Z$ where the evolution is driven by urns of
the preceding type on each integer. For other models see
\cite{Pem2}. A first way to define it is as follows. Let $f:[0,1]\to
(0,1)$ and $(\alpha_0^x,l_0^x)_{x\in \Z} \in ([0,1]\times
(0,+\infty))^\Z$ be given. Then if $x\in \Z$, set $L_0^x=0$ and for $n\ge 1$ let
$$L_n^x:=\sum_{k=0}^{n-1} 1_{\{X_k=x\}},$$
be the total time spent in $x$ up to time $n-1$ by the random walk. Let then
$$\tilde{\alpha}_n^x:= \frac{1}{l_0^x+L_{n}^x}\left\{\alpha_0^xl_0^x+\sum_{k=0}^{n-1}
1_{\{X_k=x,X_{k+1}=x+1\}}\right\},$$ be the proportion of times it
has moved to the right (up to some initial weights). Now if $X_n=x$,
for some $x\in \Z$ and $n\ge 0$, then $X_{n+1}=x+1$ with probability
$f(\tilde{\alpha}_n^x)$ and $X_{n+1}=x-1$ with probability
$1-f(\tilde{\alpha}_n^x)$. This defines recursively the random walk.
Moreover for $n\ge 1$, define $\tau_n^x$ as the time of the
$n^{\textrm{th}}$ return to $x$: for $n\ge 1$, $\tau_n^x=\inf\{k>\tau_{n-1}^x:~X_k=x\}$ with the convention $\inf\emptyset=\infty$, and $\tau_0^x=\inf\{k\ge 0:~X_k=x\}$. Set also $l_n^x:=l_0^x+n$, for $n\ge 1$. Let
$$\alpha_n^x:=\tilde{\alpha}_{\tau_{n-1}^x+1}^x,$$
when $\tau_{n-1}^x<\infty$, and $\alpha_n^x=0$ otherwise.
Then the processes $((\alpha_n^x,l_n^x),0\le n\le L_\infty^x)$  form a family of urn processes of the type described above stopped at the random time $L_\infty^x$. More precisely, $\{L_\infty^x>n\}=\{\tau_n^x<\infty\}\in\kF_{\tau_n^x}$ and on this event,
\begin{eqnarray*}
\pp\left[\alpha_{n+1}^x=\frac{l_n^x\alpha_n^x+1}{l_n^x+1}\Big|\kF_{\tau_n^x}\right]&=&\pp[X_{n+1}=x+1|\kF_{\tau_n^x}]\\
&=& f(\tilde{\alpha}_{\tau_n^x}^x)=f(\alpha_n^x).
\end{eqnarray*}
In the case the walk is recurrent, these urn processes are independent, and they are identically
distributed when $(\alpha_0^x,l_0^x)$ does not depend on $x$.

\vspace{0.2cm} \noindent There is another way to define this random
walk which goes in the other direction. Assume first that we are
given a family of independent urn processes $((\alpha^x_n,l^x_n),
n\ge 0)_{x\in \Z}$ indexed by $\Z$. One can consider this as the
full environment for instance. Then given the full environment the
random walk evolves deterministically: first it starts from $X_0$. Next
let $n\ge 0$ be given and assume that the random walk has been
defined up to time $n$. Suppose that $X_n=x$ and $L_n^x=k$, for some
$k\ge 0$. Then $X_{n+1}=x+1$ if $\alpha^x_{k+1}>\alpha^x_k$, and
$X_{n+1}=x-1$ otherwise.

\vspace{0.2cm} \noindent For $n\ge 0$, we define the environment
$w_n\in ([0,1]\times (0,+\infty))^\Z$ at step $n$ by
$$ w_n:= (\alpha_{L_n^x}^x,l^x_0+L^x_n)_{x\in \Z}.$$
We denote by $\kF_n$ the $\sigma$-algebra generated by
$(X_0,w_0,\dots,w_n)$, or equivalently by $(w_0,X_0,\dots,X_n)$. For
$x\in \Z$ and $w$ some environment, we denote by $\E_{x,w}$ the law
of the random walk starting from $X_0=x$ and with initial
environment $w_0=w$. If no ambiguity on $x$ or $w$ is possible we
will sometime forget them in the notation. A random walk of law
$\E_{x,w}$ will be called a generalized reinforced random walk
started at $(x,\omega)$ associated to $f$.

Observe that $((w_n,X_n),n \ge 0)$ is a Markov process (whereas
$(X_n,n\ge 0)$ is not), and in particular:
$$\E_{x,w} [g(X_{n+1}) \mid \kF_n] = \E_{X_n,w_n}[g(X_1)],$$
for any $(x,w)$ and any bounded measurable function $g$.

Note that the directed reinforced random walk started at $x_0$, with
initial weight $(a_0(x,y); x\in\Z, ~y\in\{x-1,x+1\})$ and reinforcement
parameter $\Delta$  defined in the introduction has law
$\E_{x_0,w_0}$ with $f(x)=x$ and $w_0^x=(\alpha_0^x,l_0^x)$ defined by
$$\alpha_0^x=\frac{a_0(x,x+1)}{a_0(x,x-1)+a_0(x,x+1)},$$
and $$l_0^x=\frac{a_0(x,x-1)+a_0(x,x+1)}{\Delta}.$$ The undirected reinforced random walk defined in the
introduction has also the law of a certain generalized reinforced random walk. For example, the undirected
reinforced random walk started at $0$ with initial weights $b_0(x,x+1)=b_0>0$ and reinforcement parameter
$\Delta$ has law $\E_{0,\omega_0^x}$, with $w_0^x=(\alpha_0^x,l_0^x)$ defined by
\begin{eqnarray*}
w_0^0=\left(\frac{1}{2},\frac{b_0}{\Delta}\right)\quad \hbox{and} \quad 
w_0^x=\left\{
\begin{array}{ll}
\left(\frac{b_0}{2b_0+\Delta},\frac{2b_0+\Delta}{2\Delta}\right) & \hbox{if } x\geq 1,\\
\left(\frac{b_0+\Delta}{2b_0+\Delta},\frac{2b_0+\Delta}{2\Delta}\right) & \hbox{if } x\leq -1.
\end{array}
\right.
\end{eqnarray*}
This corresponds to the case when $l_0^x=l_0\in
(0,+\infty)$ for all $x\neq 0$, $\alpha_0^x=\alpha_0\in (0,1)$ for
$x\geq 1$, and $\alpha_0^x=1-\alpha_0$, for $x\leq -1$.

\medskip
In the following $w_0$ will satisfy:
\begin{Hypo} \label{hypw+}
The starting environment is such that for all $x\ge 1$,
$w_0^x=w_0^1$.
\end{Hypo}
or will satisfy: 
\begin{Hypo} \label{hypw}
The starting environment is such that for all $x\ge 1$,
$w_0^x=w_0^1$ and for all $x\leq -1$, $w_0^x=w_0^{-1}$.
\end{Hypo}

\subsection{Hypothesis on $f$ and stable points}
\label{secf} Throughout the paper $f:[0,1]\to(0,1)$ will be a
regular function ($C^3$ is enough for our purpose). We say that $p$
is a fixed point if $f(p)=p$. It is called stable if $f'(p)\le 1$.
We will assume that all fixed points of $f$ are isolated.

\subsection{Statement of the main results}
Let $X$ be a reinforced random walk of law $\pp_{0,w_0}$, for some initial environment $w_0$.
This walk is called recurrent
if it visits every site infinitely often, and it is
called transient if it converges to $+\infty$ or to
$-\infty$. We denote by $R$ the event of recurrence, 
%$R=\cap_{x\in\Z} R_x$ with $R_x=\{L_\infty^x=\infty\}$, 
and by $T$ the event of transience, $T=\{X_n\to+\infty\}\cup\{X_n\to -\infty\}$. In section \ref{sec:01}, it will be shown that, under Hypothesis \ref{hypw} (or under Hypothesis \ref{hypw+} if $f\ge 1/2$), $X$ is either a.s. recurrent or a.s. transient.

The drift accumulated at time $n$ by $X$ is equal to
$$\sum_x\sum_{k=0}^{L_n^x-1} (2f(\alpha_k^x)-1).$$ 
The methods developed in this paper are well adapted to the particular case $f\geq 1/2$, making this drift nonnegative and
nondecreasing. In this case one can define for all $x\in \Z$, 
$$\delta_\infty^x:=\sum_{k=0}^\infty (2f(\alpha_k^x)-1),$$
which is the drift accumulated at site $x$ if the random
walk visits $x$ infinitely often. Then we have
\begin{theo}\label{theo1}
Assume Hypothesis \ref{hypw+} and that $f\ge 1/2$. Then the random walk $(X_n, n\ge 0)$ is
recurrent if, and only if, $\E[\delta^1_\infty]\le 1$.
\end{theo}
Note that this theorem is an analogue of Zerner's criterion \cite{Zer} for cookie
random walks.
Using the results of section \ref{securnes} where the finiteness of $\E[\delta^1_\infty]$ is discussed, this theorem implies in particular
\begin{cor} Assume Hypothesis \ref{hypw+} and that $f\ge 1/2$.
If $1/2$ is not the unique stable fixed point of $f$, or if
$f''(1/2) >0$, then $(X_n, n\ge 0)$ is a.s. transient.
\end{cor}
\begin{proof} Let $p$ be a stable fixed point of $f$, with $p\neq 1/2$.
As states Theorem \ref{convpfixes} below, $\alpha_k^1$ converges to $p$ with positive probability. Thus, with positive probability
$\delta_\infty^1=\infty$. When $p=1/2$ is the only fixed point and if $f''(1/2)>0$, Proposition \ref{driftinfini} below shows that   $\delta_\infty^1=+\infty$ a.s. and we conclude by using Theorem \ref{theo1}. 
\end{proof}

Without the assumption that $f\ge 1/2$ 
%(under the assumptions of the next theorem, $\delta^x_\infty$ and $\E[\delta^x_\infty]$ will be well defined
%\footnote{when $1/2$ is the only fixed point and when $f'(1/2)\neq 0$, $\delta^x_\infty$ is not clear defined.})
, we will prove the
\begin{theo} \label{theo2}
Assume Hypothesis \ref{hypw} and that $f$ has a unique fixed point $p$.
\begin{itemize}
\item If $p\neq 1/2$ then $\pp[R]=0$.
\item If $p=1/2$ and $f'(1/2)=0$, then $\E[\delta^1_\infty]>1$ or $\E[\delta^{-1}_\infty]<-1$ imply $\pp[R]=0$.
\end{itemize}
\end{theo}
Notice, what is part of the result, that $\delta^x_\infty$ and $\E[\delta^x_\infty]$ are still well defined for all $x$ under the hypothesis of the theorem.

The sufficient condition to get $\pp[R]=0$ in the case $p=1/2$ has
to be compared to the result of \cite{KZer} in the context of cookie
random walks, where it is proved that this is also a necessary
condition. Here we were not able to prove this.

Theorem \ref{theo2} implies in particular 
\begin{cor} Assume Hypothesis \ref{hypw} and that $1/2$ is the only fixed point of $f$. If $f'(1/2)=0$ and $f''(1/2)\neq 0$, then
$\pp[R]=0$.
\end{cor}
\begin{proof} This follows from Proposition \ref{driftinfini} which shows that $\delta_\infty^1=+\infty$ a.s. if $f''(1/2)>0$, and $\delta_\infty^{-1}=-\infty$ a.s. if $f''(1/2)<0$.
\end{proof} 

These results allow to describe interesting phase transitions. This
will be done in the last section. For example, there exists a
function  $f\ge 1/2$ having $1/2$ as a unique stable fixed point, such
that if $X$ has law $\pp=\pp_{0,w_0}$, with $w_0^x=(1/2,l)$ for all
$x$, then
\begin{theo}
\label{TC} There exists $l_1>0$ such that if $l\ge l_1$ the walk is
recurrent whereas if $l<l_1$ it is transient.
\end{theo}
This phase transition is similar, yet opposite, to the one observed by Pemantle for edge-reinforced random walks
on trees \cite{Pem000}: there exists $\Delta_1>0$ such that if the  reinforcement parameter $\Delta$ is smaller
than $\Delta_1$ it is transient whereas it is recurrent  when this parameter is greater than $\Delta_1$. Indeed, in
the non-oriented reinforced framework discussed in the introduction, starting with small $l$ is equivalent to starting with
large $\Delta$.

%%%%%%%%%%%%%%%%%%%%%%%%%%%%%%%%%%%%%%%%%%%%%%%%%%%%%%%%%%%%%%%%%%%%%%%%%%%%%%%%%%%%%%%%%%%%%%%%%%%%%%%%%%%%%%%%%%%%%

\section{Preliminaries on urns}

\subsection{Convergence of urn processes}

We recall here some known results about convergence of urn
processes. In particular the next theorem is of fundamental
importance in all this paper. Remember that all functions $f$
considered in this paper satisfy the hypothesis of section
\ref{secf}.
\begin{theo}[\cite{HLS}, \cite{Pem1}]
\label{convpfixes} Let $(\alpha_n,n\ge 0)$ be an urn process
associated to some function $f$. Then almost surely $\alpha_n$
converges to a stable fixed point of $f$ and for any stable fixed
point, the probability that $\alpha_n$ converges to this point is
positive.
\end{theo}
The convergence to a stable fixed point $p$ with positive
probability was first proved in \cite{HLS}, when $f(x)-x$ changes of
sign near $p$, and in \cite{Pem1} in the special case when the sign
of $f(x)-x$ is constant near $p$. The non existence a.s. of other
limiting points was also first proved in \cite{HLS} (for extensions
to more general settings see \cite{D}, \cite{B}, \cite{Pem00},
\cite{Pem0}).

\vspace{0.1cm} \noindent There is also a central limit theorem,
which can be extracted from the book of Duflo:

\begin{theo}[\cite{D} Theorem 4.III.5]
\label{tcl} Suppose that $p\in (0,1)$ is a stable fixed point of
$f$. Let $a=f'(p)$ and $v^2=p(1-p)$. If $a<1/2$, then conditionally
on $\alpha_n \to p$, $\sqrt{n}(\alpha_n-p)$ converges in law, as $n$
tends to $+\infty$, toward a normal variable with variance
$v^2/(1-2a)$.
\end{theo}

\subsection{Convergence of the drift}
\label{securnes}
For $n\in \N$, we set
$$\delta_n= \sum_{k=0}^n (2f(\alpha_k)-1).$$
Then $\delta_n$ will correspond to the drift accumulated at a given
site after $n+1$ visits to this site. If $\delta_n$ converges when
$n\to + \infty$, we denote by $\delta_\infty$ its limit. We will
need also to consider its negative and positive parts defined
respectively by
$$\delta^-_n:= \sum_{k=0}^n (2f(\alpha_k)-1)_-,$$
and
$$\delta^+_n:= \sum_{k=0}^n (2f(\alpha_k)-1)_+,$$
for all $n\ge 0$. In fact we can always define in the same way
$\delta_\infty^-$ and $\delta_\infty^+$, even when $\delta_n$ does
not converge. Moreover observe that if $\E[\delta_\infty^-]$ or
$\E[\delta_\infty^+]$ is finite, then $\delta_n$ converges a.s. and
$\E[\delta_\infty]=\E[\delta_\infty^+]-\E[\delta_\infty^-]$. It
happens that for our purpose, such finiteness result will be needed.

\vspace{0.2cm} The problem is that the convergence of the drift
appears to be a non-trivial question. To be more precise, we were
able to obtain a satisfying result essentially only 
when $f$ has a unique fixed point. When this fixed point is greater
(resp. smaller) than $1/2$, it is immediate to see that the drift
converges a.s. toward $+\infty$ (resp. $-\infty$). However to see
that $\E[\delta_\infty^-]$ (resp. $\E[\delta_\infty^+]$) is finite,
some non-trivial argument is needed. 
Since it is the same as in
the more difficult case when $1/2$ is the unique fixed point, we start by this case.
Let us give here an heuristic of how we handle this convergence problem when $p=1/2$:
the central limit theorem (Theorem \ref{tcl}) shows that $\sqrt{k}(\alpha_k-1/2)$ converges in law.
A Taylor expansion of $f$ shows that
\begin{equation} 
\label{taylor}
2f(\alpha_k)-1=2f'(1/2)(\alpha_k-1/2)+f''(1/2)(\alpha_k-1/2)^2+O((\alpha_k-1/2)^3).
\end{equation}
The first term is of order $k^{-1/2}$, the second of order $k^{-1}$  and the third of order $k^{-3/2}$.
When $f'(1/2)\neq 0$, then $\E[\delta_\infty^-]=\E[\delta_\infty^+]=\infty$ (see Proposition \ref{f'12} below).
When $f'(1/2)=0$ and $f''(1/2)>0$, $\delta_\infty^+=\infty$ and since $(2f(\alpha_k)-1)_-$ is of order $k^{-3/2}$, $\E[\delta_\infty^-]<\infty$. Finally, when $f'(1/2)=f''(1/2)=0$, both $\E[\delta_\infty^-]$ and $\E[\delta_\infty^+]$ are finite (see Proposition \ref{1/2pfixe} below).
\begin{prop}
\label{f'12}
Assume that $1/2$ is the unique fixed point of $f$
and that $f'(1/2)\neq 0$. Then $\E[\delta_\infty^-]=\E[\delta_\infty^+]=+\infty$.
\end{prop}
\begin{proof}
Let us prove that $\E[\delta_\infty^+]=+\infty$ (the proof of the other equality $\E[\delta_\infty^-]=+\infty$ is identical). 
Since $f'(1/2)\neq 0$, there exist positive constants $c_1$, $c_2$, such that 
\begin{eqnarray}
\E\left[(2f(\alpha_k)-1)_+\right]& \ge & c_1\E\left[(\alpha_k-1/2)_+1_{\{(\alpha_k-1/2)<c_2\}}\right] \\
                                   & \ge & \frac{c_1}{\sqrt{k}} \E\left[(\sqrt{k}(\alpha_k-1/2))_+1_{\{\sqrt{k}(\alpha_k-1/2)<1\}}\right],
\end{eqnarray}
for $k$ large enough. Then the central limit theorem (Theorem \ref{tcl}) gives $\sum_k \E[(2f(\alpha_k)-1)_+]=+\infty$. This proves that $\E[\delta_\infty^+]=+\infty$, as wanted.  
\end{proof}

\begin{prop}
\label{1/2pfixe} Assume that $1/2$ is the unique fixed point of $f$
and that $f'(1/2)=0$. Then
\begin{itemize}
\item[$\bullet$] If $f''(1/2)>0$, then $\E[\delta_\infty^-]<+\infty$.
\item[$\bullet$] If $f''(1/2)<0$, then $\E[\delta_\infty^+]<+\infty$.
\item[$\bullet$] If $f''(1/2)=0$, then $\E[\delta_\infty^-]$ and  $\E[\delta_\infty^+]$ are both finite.
\end{itemize}
In all cases, $\delta_n$ converges a.s. toward $\delta_\infty$ and
$\E[\delta_\infty]$ is well defined.
\end{prop}

\begin{proof}
For $x\in [-1/2,1/2]$, let $h(x):=f(x+1/2)-x-1/2$, and for $n\ge 0$,
let $x_n:=\alpha_n-1/2$. Let also $\epsilon_{n+1}$ be equal to $1$
if $\alpha_{n+1}>\alpha_n$ and equal to $-1$ otherwise. By
definition $(x_n,n\ge 0)$ satisfies the following stochastic
algorithm:
\begin{eqnarray}
\label{algostoc} x_{n+1}=x_n + \frac{h(x_n)}{l_0+n+1} +
\frac{\xi_{n+1}}{l_0+n+1},
\end{eqnarray}
where $\xi_{n+1}=\epsilon_{n+1}-\E[\epsilon_{n+1}\mid \kF_n]$.

\smallskip
\noindent$\bullet\quad$ Consider first the case $f''(1/2)\neq 0$. Then the sign of
$f-1/2$ is constant in a neighborhood of $1/2$. To fix ideas let say
that $f\ge 1/2$ in $[1/2-\epsilon, 1/2 +\epsilon]$ for some constant
$\epsilon >0$. In this case we will prove that
$\E[\delta_\infty^-]<+\infty$. For all $n\ge 0$,
$$ \E[\delta_n^-] \le \sum_{k=0}^n \pp[x_k^2\ge \epsilon^2].$$
Therefore it suffices to prove that this last series is convergent.
This will be achieved by using Equation \eqref{algostoc} and some
ideas from the proof of Proposition 4.1.3  in \cite{D}. First, since
$1/2$ is the unique fixed point of $f$, there exists some constant
$a>0$ such that $xh(x) \le -ax^2/2$ for all $x\in [-1/2,1/2]$.
Moreover we can always take $a<1/2$. Next Equation \eqref{algostoc}
gives
\begin{eqnarray*}
x_{n+1}^2 \le x_n^2(1-\frac{a}{l_0+n+1}) +
\frac{2x_n\xi_{n+1}}{l_0+n+1} + \frac{u_{n+1}}{(l_0+n+1)^2},
\end{eqnarray*}
where $u_{n+1}$ is bounded, i.e. there exists $C>0$ such that a.s.
$|u_{n+1}|\le C$ for all $n$. Let $s_n = \prod_{k=0}^n
(1-a/(l_0+k+1))$. By induction we get
$$x_n^2 \le s_n x_0^2 + s_n \sum_{k=0}^n ((l_0+k+1)s_k)^{-1}x_k \xi_{k+1} + C s_n \sum_{k=0}^n ((l_0+k+1)^2s_k)^{-1}.$$
Since $s_k \sim k^{-a}$, this gives
\begin{eqnarray*}
x_n^2 \le C' s_n +s_n \sum_{k=0}^n ((l_0+k+1)s_k)^{-1}x_k \xi_{k+1},
\end{eqnarray*}
for some constant $C'>0$. Define the martingale $(M_n,n\ge 0)$ by
$$M_n = \sum_{k=0}^n ((l_0+k+1)s_k)^{-1}x_k \xi_{k+1} \quad \forall n\ge 0.$$
Then for $n$ large enough and any integer $\alpha >0$, we have
$$\pp[x_n^2\ge \epsilon^2] \le \pp[s_n|M_n| \ge \epsilon^2/2]\le \frac{(2s_n)^{2\alpha}}{\epsilon^{4\alpha}}
\E[M_n^{2\alpha}].$$ But Burkholder-Davis-Gundy inequality (see e.g.
\cite{W} p.151) implies that for some constant $c_\alpha$,
$$\E[M_n^{2\alpha}]\le c_\alpha \E[<M>_n^\alpha],$$
where $<M>_n := \sum_{k=0}^n ((l_0+k+1)s_k)^{-2}x_k^2 \xi_{k+1}^2$.
Moreover, since $a<1/2$,
$$<M>_n \le \sum_{k=0}^n ((l_0+k+1)s_k)^{-2}\le C,$$
for some constant $C>0$. Thus
$$ s_n^{2\alpha}E[<M>_n^\alpha]\le C_\alpha n^{-2a\alpha},$$
with $C_\alpha>0$ a constant. Taking now $\alpha$ large enough shows
that
$$\sum_{k=0}^{+\infty} \pp[x_n^2\ge \epsilon^2]<+\infty,$$
as we wanted.

\smallskip
\noindent$\bullet\quad$ It remains to consider the case when $f''(1/2)=0$. There exists a constant $C>0$ such that
$$|2f(\alpha_n)-1| \le C |x_n|^3 \quad \forall n\ge 0,$$
from which we get
$$\delta_n^-+\delta_n^+ \le C \sum_{k=0}^n |x_k|^3.$$
Thus it suffices to prove that $\sum_{k=0}^{+\infty}\E[|x_k|^3]$ is
finite. But since $f'(1/2)=0$, there exists $\epsilon>0$ such that
$2xh(x) \le -x^2$ when $|x|\le \epsilon$. Therefore \eqref{algostoc}
gives in fact
$$\E[x_{n+1}^2] \le \E[x_n^2](1-\frac{1}{n}) + 4\frac{\pp[|x_n|\ge \epsilon]}{n} + \frac{C}{n^2}.$$
Now the proof of the preceding case shows that
$$\E[x_{n+1}^2] \le \E[x_n^2](1-\frac{1}{n}) + \frac{C'}{n^2}.$$
This proves by induction that $\E[x_n^2]\le C/n$ for some constant
$C>0$. Let us now consider the moments of order $4$. Since $4x^3h(x)
\le-3x^4$ in $[-\epsilon, \epsilon]$ for some $\epsilon>0$,
\eqref{algostoc} gives similarly
$$\E[x_{n+1}^4] \le \E[x_n^4](1-\frac{3}{n}) + \frac{C}{n^2},$$
for some constant $C>0$. By induction, this gives $E[x_n^4] \le C'
n^{-2}$, with $C'>0$ another constant. Then Cauchy-Schwarz
inequality gives (up to constants)
$$\E[|x_n|^3]\le (\E[x_n^2]\E[x_n^4])^{1/2}\le n^{-3/2},$$
which is summable. This finishes the proof of the proposition.
\end{proof}

Observe that the argument given in the proof of the above
proposition applies as well when the unique fixed point of $f$ is
different from $1/2$. Thus we proved also the

\begin{prop}
\label{pfixeautre} If $f$ has a unique fixed point $p>1/2$,
resp. $p<1/2$, then $\E[\delta_\infty^-]$, resp.
$\E[\delta_\infty^+]$, is finite. In particular $\delta_n$ converges
a.s. toward $+\infty$, resp. $-\infty$.
\end{prop}

Our last result concerns the a.s. non-finiteness of $\delta_\infty$.
First if $p\neq 1/2$ is a stable fixed point of $f$, then
conditionally on $\{\alpha_n \to p\}$, $\delta_n/n$ converges toward
$2p-1$ and thus $\delta_\infty=+\infty$. The next result
investigates the case $p=1/2$.
\begin{prop}
\label{driftinfini} If $1/2$ is a stable fixed point of $f$,
$f'(1/2)=0$ and $f''(1/2)>0$ (respectively $f''(1/2)<0$), then
conditionally on $\alpha_n \to 1/2$, almost surely
$\delta_\infty=+\infty$ (respectively $\delta_\infty=-\infty$).
\end{prop}
\begin{proof}
To fix ideas assume that $f''(1/2)>0$. The other case is analogous. A
limited development of $f$ near $1/2$ gives
\begin{eqnarray}
\label{delta} \left|\delta_n - f''(1/2) \sum_{k=0}^n (\alpha_k-1/2)^2\right|
\le C \sum_{k=0}^n |\alpha_k-1/2|^3 \quad \forall n\ge 0,
\end{eqnarray}
with $C>0$ some positive constant. For $n\ge 0$, we set $z_n:=
\sqrt{n}(\alpha_n-1/2)$. We already saw in Theorem \ref{tcl} that conditionally on $\{\alpha_n \to 1/2\}$, 
$z_n$ converges in law toward a normal variable. In fact this holds
in the sense of the trajectory. More precisely, an elementary
calculus shows that $(z_n,n\ge 0)$ is solution of a stochastic
algorithm of the form:
\begin{equation}\label{eq:z}z_{n+1}=z_n-\frac{z_n/2 - r_{n+1}}{l_0+n+1} + \frac{\xi_{n+1}}{\sqrt{l_0+n+1}},\end{equation}
where $r_{n+1}=\kO(\sqrt{n}(\alpha_n-1/2)^2+n^{-1})$. For $t\in
[\log n, \log (n+1)]$, let
\begin{equation}\label{eq:Y}Y_t= z_{n+1} +(t-\log n) (-z_{n+1}/2+r_{n+1})+ (t-\log n)^{1/2}\xi_{n+2}.\end{equation}
For $u\ge 0$, call $(Y^{(u)}_t,t\ge 0)$ the continuous time process
defined by $Y^{(u)}_t=Y_{u+t}$ for $t\ge 0$. Then Theorem 4.II.4 in
\cite{D} says that (conditionally on $\{\alpha_n \to 1/2\}$) the sequence of processes $(Y^{(u)}_t,t\ge 0)$ converges in
law in the path space toward an Ornstein-Uhlenbeck process
$(U_s,s\ge 0)$, when $u\to +\infty$ (the condition on $r_n$ in the
hypothesis of the theorem is not needed here, as one can see with
Theorem 4.III.5 and its proof in \cite{D}). Now we will deduce from
this result that a.s. on the event $\{\alpha_n \to 1/2\}$, 
\begin{equation}
\label{eq:serie}
\sum_{k=1}^\infty \frac{z_k^2}{k} = + \infty.
\end{equation}
If we define $z_t$ for all $t\ge 0$ by $z_t=z_{[t]}$, then one can check
that $\sum_{k=0}^\infty z_k^2/k$ is comparable with
$\int_{1}^{+\infty} z_t^2/t \ dt=\int_1^\infty z_{e^t}^2 \ dt$.
So if this series is finite, then $\int_n^{n+1} z_{e^t}^2\ dt \to 0$ when $n\to +\infty$.
Moreover, using (\ref{eq:z}) and (\ref{eq:Y}), we have that a.s. on the event $\{\alpha_n\to 1/2\}$,  
$Y^2_t=z^2_{e^t}+o(1)$. Therefore a.s. on the event $\{\sum_k z_k^2/k<\infty\}\cap \{\alpha_n\to 1/2\}$, we have $\int_n^{n+1} Y_t^2\ dt \to 0$ when $n\to +\infty$. But
this cannot hold since on the event $\{\alpha_n\to 1/2\}$, $\int_n^{n+1} Y_t^2\ dt$ converges in
law toward $\int_0^1 U_s^2\ ds$, which is a.s. non-zero. Thus \eqref{eq:serie} holds.  
Finally, using the fact that a.s. on the event $\{\alpha_n\to 1/2\}$, $\sum_{k=0}^n|\alpha_k-1/2|^3=o\left(\sum_{k=1}^n \frac{z_k^2}{k}\right)$, 
\eqref{delta} and \eqref{eq:serie} show that $\delta_\infty=+\infty$.  
\end{proof}

%%%%%%%%%%%%%%%%%%%%%%%%%%%%%%%%%%%%%%%%%%%%%%%%%%%%%%%%%%%%%%%%%%%%%%%%%%%%%%%%%%%%%%%%%%%%%%%%%%%%%%%%%%
%%%%%%%%%%%%%%%%%%%%%%%%%%%%%%%%%%%%%%%%%%%%%%%%%%%%%%%%%%%%%%%%%%%%%%%%%%%%%%%%%%%%%%%%%%%%%%%%%%%%%%%%%%

\section{A zero-one law}\label{sec:01}

In all this section $X$ is a generalized reinforced random
walk of law $\pp=\pp_{0,w_0}$ associated to $f$, and $w_0$ satisfies
Hypothesis \ref{hypw}. We will try to relate its asymptotic
behavior with urn characteristics. Our first result is general. It
is a zero-one law for the property of recurrence. Remember that the
random walk is recurrent if all sites are visited infinitely often.

\begin{lem}
\label{01} Assume Hypothesis \ref{hypw}. Then $\pp[R] \in
\{0,1\}$, and $\pp[T]=1-\pp[R]$.
\end{lem}
\begin{proof}
First Borel-Cantelli Lemma implies that if a site is visited
infinitely often, then the same holds for all sites. So there are
only three alternatives. Either the random walk is recurrent, or it
tends toward $+\infty$ or toward $-\infty$. In other words $\pp[T]=1-\pp[R]$. Thus, if
$$T_n=\inf\{ k\ge 0 \mid X_k = n\} \quad \forall n\ge 0,$$
then $1_{\{T_n<+\infty\}}$ converges toward $1_{R\cup \{X_n \to
+\infty\}}$, when $n\to +\infty$. In the same way the event $\{X_n
>0 \quad \forall n>0\}$ is included in $\{X_n \to + \infty\}$. In
fact there is a stronger relation:

\begin{lem}
\label{transient} For any initial environment $w_0$ and any $k\ge 0$, $\pp[X_n\to +\infty]>0$ if, and
only if, $\pp_{k,w_0}[X_n >k; \quad \forall n> 0]
>0$.
\end{lem}

\begin{proof} We do the proof for $k=0$. The other cases are identical. This proof is similar to Zerner's proof of Lemma $8$ in \cite{Zer}. We just have to prove the only if
part. Call $\tau_2$ the last time the random walk visits the integer
$2$. If $\kC$ is some path of length $k$ starting from $0$ and
ending in $2$ on $\Z$, call $E_\kC$ the event that the random walk
follows the path $\kC$ during the first $k$ steps. Define also
$w_{\kC}$ as the state of all urns once the walker has performed the
path $\kC$. If $\pp[X_n\to +\infty]>0$, then for some path $\kC$
from $0$ to $2$, we have
$$0<\pp[E_\kC, X_n >2 \quad \forall n> k]=\pp[E_\kC]\times \pp_{2,w_\kC}[X_n>2 \quad \forall n>0].$$
Now construct $\kC'$ as follows: it starts by a jump from $0$ to $1$
and then we add (in chronological order) all the excursions of $\kC$
above level $1$. Then clearly
$$ \pp_{2,w_\kC}[X_n>2 \quad \forall n>0]=\pp_{2,w_{\kC'}}[X_n>2 \quad \forall n>0].$$
Moreover, since the range of $f$ is in $(0,1)$, it is elementary to
see that $\pp[E_{\kC'}]>0$. Thus
$$\pp[X_n >0 \quad \forall n>0] \ge  \pp[E_{\kC'}] \times \pp_{2,w_{\kC'}}[X_n>2 \quad \forall n>0] >0,$$
which proves the lemma.
\end{proof}

\noindent We can finish now the proof of Lemma \ref{01}. The
martingale convergence theorem and the Markov property imply
\begin{eqnarray*}
1_{\{X_n \to + \infty\}} & = & \lim_{n\to +\infty} \pp_{0,w_0}[X_n \to + \infty \mid \kF_{T_n}]1_{\{T_n< +\infty\}} \\
             & = & \lim_{n\to +\infty} \pp_{n,w_{T_n}} [X_n \to + \infty ]1_{\{T_n< +\infty\}}\\
            &\ge & \limsup_{n\to +\infty}\pp_{n,w_{T_n}} [X_m>n \quad \forall m > 0]1_{\{T_n< +\infty\}}\\
            & = &  \pp_{1,w_0} [X_n>1 \quad \forall n > 0]1_{R\cup \{X_n \to +\infty\}}.
\end{eqnarray*}
Then multiply the left and right part of this inequality by $1_R$
and take expectation. This gives
\begin{eqnarray}
\label{pp1}
\pp_{1,w_0}[X_n>1 \quad \forall n > 0] \ \pp[R] = 0.
\end{eqnarray} 
In the same way we have
$$\pp_{-1,w_0}[X_n< -1 \quad \forall n > 0] \ \pp[R] = 0.$$
These two equalities and Lemma \ref{transient} prove the lemma.
\end{proof}

\begin{rem} Let $T_0$ be the first return time to $0$. Then the usual equivalence $T_0<+\infty$ a.s. if and only if $0$
is a.s. visited infinitely often, is true. Indeed if $T_0<+\infty$
a.s. then by Lemma \ref{transient}, a.s. $X_n$ does not converge
toward $\pm \infty$. Then the $0-1$ law says that $R$ holds.
\end{rem}

\section{The case with only non-negative drift}

Here we assume that $f\ge 1/2$ and that $w_0$ satisfies Hypothesis
\ref{hypw+}. In
the following $X$ is a reinforced walk of law $\pp=\pp_{0,w_0}$. In
this case we have a more precise zero-one law.

\begin{lem}
\label{01bis} Assume Hypothesis \ref{hypw+} and $f\ge 1/2$. We have the alternative: either $(X_n,n\ge 0)$ is
almost surely transient toward $+\infty$, or it is almost surely
recurrent.
\end{lem}
\begin{proof}
Since $f\ge 1/2$, at each step the random walk has probability at
least $1/2$ to jump to the right. Thus an elementary coupling
argument (with the usual simple random walk on $\Z$) shows that a.s.
the random walk does not converge toward $-\infty$. We conclude with
\eqref{pp1} (which holds when assuming only Hypothesis \ref{hypw+}) and Lemma \ref{transient}.
\end{proof}

\begin{rem} 
\emph{We notice here that the hypothesis $f<1$ made in section \ref{secf} is not needed when $f\ge 1/2$. Indeed the only place where it is used is in the proof of Lemma \ref{transient} to show that $\pp[E_{\kC'}]>0$, but the reader can check that this is not needed when $f\ge 1/2$. This remark will be of interest for the last section.}
\end{rem}

\begin{proof}[Proof of Theorem \ref{theo1}:]
We follow essentially the proof of Theorem 12 in \cite{Zer}. Let us recall the main
lines. First we introduce some notation. For $n\ge 0$, let
$$U_n= \sum_{k\le n-1} 1_{\{X_k=0,\ X_{k+1}=-1\}},$$
and let
$$X_n^+=\sum_{k\le n-1} (X_{k+1}-X_k)1_{\{X_k \ge 0\}}. $$
A straightforward computation gives the equation
\begin{equation} X_n^+ = \max(X_n,0) - U_n \quad \forall n. \label{eqX+}\end{equation}
We define the drift $D^x_n$ accumulated in $x$ up to time $n$ by
$$D_n^x= \sum_{k=0}^{L^x_n} (2f(\alpha^x_k)-1),$$
and the drift $D^+_n$ accumulated in the non-negative integers by
$$D_n^+=\sum_{x\ge 0} D_n^x.$$
Let $(M_n^+,n\ge 0)$ be the process defined by
$$M_n^+=X_n^+ - D_n^+ \quad \forall n.$$
It is a basic fact that $(M^+_n, n\ge 0)$ is a martingale. In
particular for all $a\ge 0$ and all $n\ge 0$, using (\ref{eqX+})
with the martingale property,
$$\E[\max(X_{T_a \wedge n},0)] = \E[U_{T_a \wedge n}] + \E[D^+_{T_a \wedge n}],$$
where
$$T_a=\inf \{k \ge 0 \mid  X_k=a\}.$$
Now Lemma \ref{01bis} implies that $T_a$ is a.s. finite. Moreover
$(U_n,n\ge 0)$ and $(D^+_n,n\ge 0)$ are non-decreasing processes.
Thus letting $n$ go to $+\infty$ gives with the monotone convergence
theorem
\begin{eqnarray}
\label{equationdrift} a= \E[U_{T_a}]+ \E[D^+_{T_a}] \quad \forall a
\ge 0.
\end{eqnarray}
Moreover the Markov property shows that for any integer $x\in
[1,a]$,
$$\E[D^x_{T_a}]=\E[\E_{x,w_{T_x}}[D^x_{T_a}]]=E_{1,w_0}[D^1_{T_{a-x+1}}],$$
where the last equality holds because for all $y\ge x$,
$w^y_{T_x}=w_0^1$. Moreover $\E[D^0_{T_a}]$ and
$E_{1,w_0}[D^1_{T_{a-1}}]$ differ at most by $\E[N]$, where $N$ is
the number of visits to $0$ before the first visit to $1$. Since the
probability to jump from $0$ to $1$ is bounded away from $0$,
$\E[N]$ is finite. Therefore
\begin{eqnarray*}
\lim_{a\to +\infty} \frac{1}{a}\E[D^+_{T_a}] &=& \lim_{a\to
+\infty}\frac{1}{a}\left(\E_{0,w_0}[D^0_{T_a}]+\sum_{x=1}^a\E_{1,w_0}[D^1_{T_x}]\right)\\
&=& \E_{1,w_0}[D^1_\infty].
\end{eqnarray*}
Then \eqref{equationdrift} gives the inequality
$$\E[D^1_\infty] \le 1.$$
So if the random walk is recurrent, almost surely $D^1_\infty =
\delta^1_\infty$, and $\E[\delta^1_\infty] \le 1$. This gives
already the only if part of the theorem.

Assume now that the random walk is transient. Then Lemma
\ref{transient} shows that
$$\E[\delta^1_\infty-D_\infty^1] \ge c \E[\delta^1_\infty-\delta^1_0],$$
where $c= \pp_{1,w_0}[X_n >1 \quad \forall n> 0] >0$. Now since the
sequence $(\delta^1_n)_{n\ge 0}$ is non-decreasing, if
$\E[\delta^1_\infty-\delta^1_0]$
 was equal to $0$, this would mean that a.s. $\delta^1_n=\delta^1_0$ for all $n$. In other words the walk would evolve
 like the simple random walk,
which is recurrent. This is absurd. Thus $\E[\delta^1_\infty]>
\E[D^1_\infty]$. It remains to prove that $\E[D_\infty^1]=1$. From
\eqref{equationdrift} we see that it is equivalent to prove the
\begin{lem}
If the random walk is a.s. transient, then $\lim_{a\to + \infty}
\E[U_{T_a}]/a=0$.
\end{lem}
\begin{proof} This lemma can be proved by following the argument of Lemma $6$ in \cite{Zer}, that we reproduce here.
For $i\ge 1$, let
$$\sigma_i=\inf\{j\ge T_i\mid X_j=0\}.$$
We have $\E[U_{T_a}]=\sum_{i=0}^{a-1} \E[U_{T_{i+1}}-U_{T_i}]$. Next
$U_{T_{i+1}}-U_{T_i}\neq 0$ only on the set $A_i:=\{\sigma_i <
T_{i+1}\}$. Moreover \eqref{equationdrift} holds for any starting
environment. Thus
$$\E[U_{T_{i+1}}-U_{T_i}]\le \E[1_{A_i}\E_{0,w_{\sigma_i}}[U_{T_{i+1}}]]\le (i+1) \pp[A_i],$$
for all $i$. It remains to prove that
\begin{eqnarray}
\label{but} \frac{1}{a} \sum_{i=1}^a i \pp[A_i] \to 0,
\end{eqnarray}
when $a\to +\infty$. Let $Y_i=\pp[A_i \mid \kF_{T_i}]$. Since the
random walk is transient, the conditional Borel-Cantelli lemma
implies\footnote{since we were not able to find a reference we give
here a short proof: let $(H_n,n\ge 0)$ be the $(\kF_{T_n})_{n\ge 0}$
martingale defined by $H_n:=\sum_{i=1}^n1_{A_i}-Y_i$, for $n\ge 0$.
Let $l\ge 1$ and let $T'_l=\inf\{k\mid H_k\ge l\}$. Then $H_{n\wedge
T'_l}$ a.s. converges toward some limiting value $\alpha_l\in \R$,
when $n\to +\infty$. If a.s. only a finite number of $A_i$'s occur,
then a.s. $T'_l$ is infinite for some $l\ge 1$. This implies the
desired result.}
\begin{eqnarray}
\label{aide} \sum_{i\ge 1} Y_i < +\infty \quad \textrm{a.s.}
\end{eqnarray}
Moreover a coupling argument with the simple random walk and
standard results for this random walk show that a.s., $Y_i \le 1/i$
for all $i$. Let $\epsilon>0$. For all $i$,
$$\pp[A_i]=\E[Y_i1_{\{Y_i< \epsilon/i\}}] + \E[Y_i1_{\{Y_i\ge \epsilon/i\}}]\le \frac{\epsilon}{i} + \frac{\pp[Y_i\ge
\epsilon/i]}{i}.$$ So we can divide the sum in \eqref{but} in two
parts. One is lower than $\epsilon$ and the other one is equal to
$$\frac{1}{a}\E\left[\sum_{i=1}^a1_{\{Y_i\ge \epsilon/i\}}\right].$$
But since \eqref{aide} holds, a.s. the density of the $i\le a$ such
that $Y_i\ge \epsilon/i$ tends to $0$ when $a$ tends to $+\infty$.
Thus the preceding sum converges to $0$. This concludes the proof of
the lemma.
\end{proof}
This completes the proof of Theorem \ref{theo1}. \end{proof}

\medskip
In section \ref{ex} we will see different examples of functions
$f\ge 1/2$, symmetric with respect to $1/2$ which show in particular
that in the case when $1/2$ is the only stable fixed point and
$f''(1/2)=0$, both regimes (recurrence and transience) may appear.

\section{The case with a unique fixed point}
Here we do not assume anymore that $f\ge 1/2$, but we assume that
$f$ has a unique fixed point. The initial environment satisfies
Hypothesis \ref{hypw} and still $\pp=\pp_{0,w_0}$. 
%We prove now Theorem \ref{theo2}.

\begin{proof}[Proof of Theorem \ref{theo2}:]
The idea of the proof is the same as for Theorem
\ref{theo1}. However a priori we have to be careful when
taking limits since the drift $(D_n^+)_{n\ge 0}$ is not anymore a
non-decreasing function. But for any integer $x\ge 0$, Proposition
\ref{1/2pfixe} and Proposition \ref{pfixeautre} show that
$\E[D_n^x]$ converges toward $\E[D_\infty^x]$. In fact since
$\E[\delta_\infty^-]<+\infty$ or $\E[\delta_\infty^+]<+\infty$, if
we replace $n$ by any increasing sequence of stopping times $\tau_n$
converging toward $\tau_\infty$, these propositions show that
$\E[D_{\tau_n}^x]$ converges toward $\E[D_{\tau_\infty}^x]$. So in
fact we get
$$\lim_{n\to +\infty} \E[D_{T_a\wedge n}^+]=\E[D_{T_a}^+].$$
Next observe that for all $a\ge 0$ and $n\ge 0$, $X^+_{T_a\wedge n}
\le a$. Thus, using that $(M^+_n)_{n\ge 0}$ is a martingale, we have
$$\E[D_{T_a}^+]\le a \quad \forall a\ge 0.$$

Assume $\pp(R)=1$. Then the Markov property implies that if $1\le
x\le a$,
\begin{eqnarray}
\label{majorationdrift}
\E[D^x_{T_a}]=\E_{x,w_0}[D^x_{T_{a}}]=\E_{1,w_0}[D^1_{T_{a-x+1}}].
\end{eqnarray}
Letting $a$ tend to $+\infty$ in \eqref{majorationdrift}, and using
the fact that $D^+_{T_a}=\sum_{x=0}^a D^x_{T_a}$ gives
$$\E_{1,w_0}[D^1_\infty] \le 1.$$
Since $\pp[R]=1$, a.s. $D_\infty^1=\delta^1_\infty$ and we have
$\E[\delta^1_\infty]\le 1$. The other inequality
$\E[\delta^{-1}_\infty] \ge -1$ is similar. 
\end{proof}

Let us state now the following standard monotonicity argument:

\begin{lem}
\label{monotone} Let $f \le g$ be two functions. Then there exists a
coupling of two urn processes $((\alpha_n,l_n),n\ge 0)$ and
$((\beta_n,l'_n),n\ge 0)$ associated respectively to $f$ and $g$,
such that $l_0=l'_0$, $\alpha_0=\beta_0$, and almost surely
$\alpha_n \le \beta_n$ for all $n\ge 0$.
\end{lem}

\begin{proof} The proof is standard. Let $(U_i,i\ge 0)$ be a sequence of i.i.d. random variables uniformly distributed on
$[0,1]$. We define two urn processes starting with initial
conditions like in the lemma. Then at step $n$,
$\alpha_{n+1}>\alpha_n$ if, and only if $f(\alpha_n)\ge U_n$. The
same for $\beta_{n+1}$ (with $g$ in place of $f$). Assume now that
for some $n$, $\alpha_n>\beta_n$. Assume also that $n$ is the lowest
index where such inequality occurs. This means that
$\alpha_{n-1}=\beta_{n-1}$. But since $f\le g$, by definition of our
processes, we get an absurdity.
\end{proof}

This lemma together with Theorem \ref{theo2} allows to
consider also the case when $f$ has possibly more than one fixed
point but under the condition $f\ge 1/2$ on $[1/2,1]$. More
precisely we have

\begin{cor}
Assume that Hypothesis \ref{hypw} holds, that $f\ge 1/2$ on $[1/2,1]$ and that all fixed points of $f$
are greater or equal to $1/2$.
\begin{itemize}
\item If $1/2$ is not a fixed point, then $\pp[R]=0$.
\item If $1/2$ is a fixed point, but not the only fixed point, and $f'(1/2)=0$, then $\pp[R]=0$.
\end{itemize}
\end{cor}

\begin{proof}
If any of the two hypothesis of the corollary is satisfied, then
there exists a function $g$ such that $g\le f$, $g$ has a unique
fixed point equal to $1/2$, and $g'(1/2)=0$. We can also assume that
$g$ is increasing on $[0,1/2]$. Applying Lemma \ref{monotone} we see
that there exists an urn process $(\beta_n,n\ge 0)$ associated to
$g$ such that $\beta_n \le \alpha_n$ for all $n$. Now the proof of
Proposition \ref{1/2pfixe} shows that
$$\sum_{n\ge 0} \E[(2g(\beta_n)-1)_-] <+\infty.$$
Since $g$ is increasing on $[0,1/2]$ and $f\ge 1/2$ on $[1/2,1]$,
this implies that $\E[\delta_\infty^-]$ is finite. Moreover we know
that $\delta_\infty = +\infty$ a.s. So we have everything to apply
the proof of Theorem \ref{theo2} and to conclude.
\end{proof}

\section{Some examples}
\label{ex} Our goal here is to give examples of functions $f$
leading to interesting behavior for the associated random walk, in view of the
previous results. In all this section we consider a function $f$, symmetric
with respect to $1/2$, i.e. such that $f(1/2-x)=f(1/2+x)$ for all
$x\in [0,1/2]$, decreasing on $[0,1/2]$ and 
increasing on $[1/2,1]$. We assume also that $f$ has a unique fixed
point, equal to $1/2$, and that $f''(1/2)=0$.

\vspace{0.2cm} \noindent We start now by a comparison result. Let
$u$ be some positive real number. Define $f_u$ by the equation
$2f_u-1=\big(u(2f-1)\big)\wedge 1$. One can see immediately that
$f_u$ has the same properties as $f$ for all $u$, and moreover
that $f_u\le f_v$ if $u\le v$. Denote by $((\alpha^u_n,l_n^u),n\ge
0)$ an urn process associated to $f_u$ such that
$(\alpha_0^u,l_0^u)=(1/2,l)$ with $l>0$, and set
$\delta_\infty^u:=\sum_{n\ge 0} (2f(\alpha^u_n)-1)$. Then we have the
\begin{lem}
\label{deltau} For all $u$, $\E[\delta^u_\infty]<+\infty$. The maps
$u\mapsto \E[\delta^u_\infty]/u$ and $u\mapsto \E[\delta^u_\infty]$,
are nondecreasing respectively on $(0,1]$ and on $[0,+\infty)$. In particular $\E[\delta^u_\infty]\to 0$, when $u\to 0$.  
Moreover $\E[\delta^u_\infty]\to +\infty$, when $u\to +\infty$.
\end{lem}
\begin{proof} The first claim results from the proof of Proposition \ref{1/2pfixe}.
For the second claim, consider first $0<u<v\le 1$. By symmetry, for
any $k\ge 0$,
$$\E[2f_u(\alpha_k^u)-1]=2 \E[(2f_u(\alpha_k^u)-1)1_{\{\alpha^u_k\ge 1/2\}}].$$
Moreover, since $f_v$ is nondecreasing on $[1/2,1]$ and since one may
couple $\alpha^u_k$ and $\alpha^v_k$ such that $\alpha^u_k\le
\alpha^v_k$ a.s. by Lemma \ref{monotone},
\begin{eqnarray*}
\E[(2f_u(\alpha_k^u)-1)1_{\{\alpha^u_k\ge 1/2\}}] & = & \frac{u}{v} \E[(2f_v(\alpha_k^u)-1)1_{\{\alpha^u_k\ge 1/2\}}] \\
                        &\le & \frac{u}{v} \E[(2f_v(\alpha_k^v)-1)1_{\{\alpha^v_k\ge 1/2\}}].
\end{eqnarray*}
The result follows by summation.

The fact that $u\mapsto \E[\delta^u_\infty]$ is nondecreasing on
$[0,+\infty[$ is similar. It remains to find the limit when $u\to
+\infty$. For this, fix some $n\ge 1$. Then one can observe that
there exists $\epsilon>0$, such that $|\alpha_{2k+1}^u-1/2|\ge
\epsilon$ for all $k\le n$. This implies that for $u$ large enough,
$\E[\delta_\infty^u] \ge n/2$. Since this holds for all $n$, the
result follows.
\end{proof}
The preceding lemma and Theorem \ref{theo1} show that
there is a phase transition: let $X$ be a generalized random walk
started at $(0,w_0)$ associated to $f_u$, where the initial
environment is such that $w_0^x=(1/2,l)$. Then there exists some
$u_0>0$ such that for $u>u_0$, the random walk associated to $f_u$
is transient, whereas for $u<u_0$ it is recurrent. In particular
recurrence and transience may both appear. The question of what
happens at $u_0$ is related to the continuity of $\E[\delta_\infty]$
with respect to $f$. But explicit calculus show that
if $u\to u_0$, then for all $n$, $\E[\delta_n^u] \to
\E[\delta_n^{u_0}]$. Together with the monotonicity of
$\E[\delta_\infty^u]$ in $u$, this proves that $\E[\delta_\infty^u]$
is continuous in $u$. In particular for $u=u_0$ the random walk is
recurrent.

\vspace{0.2cm} \noindent Our second problem concerns what happens
when initial conditions vary. Here also we will see that there is
possibly a phase transition. For $(\alpha,l)\in [0,1]\times
]0,+\infty[$, we denote by $\E_{\alpha,l}$ the law of an urn process
starting from $(\alpha,l)$. First let us prove the
\begin{lem} \label{alpha} Let $\alpha\in [0,1]\smallsetminus \{1/2\}$ and $l\in (0,+\infty)$, be such that $2\alpha l-l
\in \N$. Then $\E_{\alpha,2l}[\delta_\infty] >
\E_{1/2,2l}[\delta_\infty]$.
\end{lem}
\begin{proof} We use a standard coupling argument. Let $(U_i)_{i\ge 0}$ be a family of i.i.d. random variables, uniformly
distributed in $[0,1]$. Let start two urn processes $(\alpha_n,n\ge
0)$ and $(\beta_n,n\ge 0)$, respectively from  $(\alpha,2l)$ and
$(1/2,2l)$. They evolve according to the following rule. If at step
$n$, $\alpha_n$ or $\beta_n$ is equal to $x\ge 1/2$, then we add one
Red ball in the corresponding urn if $U_n\le f(x)$. Now if $x<1/2$,
then we add a Red ball if $U_n\ge 1-f(x)$. The condition $2\alpha
l-l \in \N$ assures by induction that $l_n \alpha_n-l_n\beta_n \in
\Z$ for all $n\ge 1$. This in turn shows that the two urn processes
(as well as their symmetric with respect to $1/2$) cannot cross each
other without meeting them. Thus for all $n\ge 0$, $|\beta_n-1/2|\le
|\alpha_n-1/2|$. The lemma follows.
\end{proof}
The preceding results show in particular that the property of
recurrence or transience may depend on the initial conditions of the
urns (even if $l_0$ is fixed). Indeed it suffices to consider $f$
such that $\E_{1/2,2l_0}[\delta_\infty]=1$, which is possible by
Lemma \ref{deltau} and the continuity in $u$ of
$\E[\delta^u_\infty]$ as explained above. Then the preceding lemma
shows that for any $\alpha \neq 1/2$ satisfying the condition of the
lemma, the random walk associated with urns starting from
$(\alpha,2l_0)$ is always transient, whereas it is recurrent if they
start from $(1/2,2l_0)$.

\vspace{0.2cm} \noindent We arrive now to our last result.
\begin{lem} The map $l\mapsto \E_{1/2,2l}[\delta_\infty]$ is continuous on $(0,+\infty)$, non-increasing, and converges
toward $0$ when $l \to +\infty$.
\end{lem}
\begin{proof}
The continuity of the map is similar to the continuity of
$\E[\delta_\infty^u]$ in $u$, observed above. The fact that the map
is nonincreasing can be seen by using a coupling argument like in the
preceding lemma. Indeed let $l_0<l_1$. Let $(\alpha_n,n\ge 0)$ and
$(\beta_n,n\ge 0)$ be two urn processes starting respectively from
$(1/2,2l_1)$ and $(1/2,2l_0)$. Define their joint law like in the
previous lemma. Observe that each urn process cannot jump above
$1/2$ without touching it. In the same way, if for some $n$, $1/2\le
\beta_{n+1}<\beta_n$ and $1/2\le \alpha_{n+1}<\alpha_n<\beta_n$,
then
$$\beta_n- \beta_{n+1}=\frac{\beta_n}{2l_0+n+1},$$
and
$$\alpha_n-\alpha_{n+1}=\frac{\alpha_{n+1}}{2l_1+n}.$$
Thus
$$\beta_n-\beta_{n+1}\le \beta_n-\alpha_{n+1}.$$
In other words the two urn processes cannot cross each other without
meeting them. Thus for all $n\ge 0$, $|\beta_n-1/2|\ge
|\alpha_n-1/2|$, which proves the desired result. It remains to find
the limit when $l\to +\infty$. But for each $n$,
$\E_{1/2,2l}[\delta_n]$ converges to $0$ when $l\to +\infty$. Since
moreover for fixed $l$ it converges toward
$\E_{1/2,2l}[\delta_\infty]$, when $n\to +\infty$, the result
follows. This finishes the proof of the lemma. \end{proof}

We finish by the

\begin{proof}[Proof of Theorem \ref{TC}:] 
It suffices to choose $f$ and $l_0$ such that
$\E_{1/2,2l_0}[\delta_\infty]>1$. Then the result follows
immediately from the preceding lemma.
\end{proof}

\end{document}